
\input amstex
\documentstyle{amsppt}

\NoBlackBoxes
\topmatter
\define\osr{\overset \sim\to\rightarrow}

\define \Dl{\Delta}

\define \G{\Gamma}

\define \vp{\varphi}

\define \iy{\infty}
\define \la{\langle}
\define \ra{\rangle}

\define \tm{\times}

\define \ri{\rightarrow}

\define \sbt{\subset}

\define \sbq{\subseteq}

\define \ov{\overline}

\define \edm{\enddemo}
\define \ep{\endproclaim}

\define \1{^{-1}}
\define \2{^{-2}}

\define \CP{\Bbb C\Bbb P}

\define \CPt{\Bbb C\Bbb P^2}
\define \BC{\Bbb C}

\define \BZ{\Bbb Z}

\define \tB{\tilde{B}}

\define \tT{\tilde{T}}
\define \tX{\tilde{X}}

\define \tP{\tilde{P}}

\define \BCen{\operatorname{Center}}

\define \Xab{X_{ab}}
\baselineskip 20pt

\topmatter
\title{Fundamental groups of complements of\\ branch curves as solvable
groups}\endtitle
\author{ Boris Moishezon
\address Department of Mathematics,Columbia University, New York, NY
\endaddress
 and Mina Teicher}\endauthor
\address Deptartment of Mathematics and Computer Science,
Bar-Ilan University, Ramat-Gan  52900, Israel\endaddress
\abstract
In this paper we show that fundamental groups of complements of curves are
``small'' in the sense that they are ``almost solvable''.  Thus we can start to
compute $\pi_2$ as a module over $\pi_1$ in order to produce new invariants of
surfaces that might distinguish different components of a moduli space.
\endabstract \subjclass 14J10, 20F36\endsubjclass \keywords Fundamental groups,
branch curves, Veronese surfaces\endkeywords \leftheadtext{B. Moishezon and M.
Teicher} \rightheadtext{Fundamental groups of branch curves}

\endtopmatter
\document
\subhead{0.  Applications of the calculations of fundamental groups to
algebraic surfaces}\endsubhead

Our study of fundamental groups of branch curves is aimed towards understanding
algebraic surfaces.

Algebraic surfaces are classified by discrete and continuous invariants.
Fixing the discrete invariants (of homotopy type), one gets a family of
algebraic
surfaces parametrized by an algebraic variety which is called the moduli space
(the word ``moduli" stands for ``continuous invariants").
Very little is known about the structure of moduli
space of surfaces of general type.

	We study algebraic surfaces in order to understand the
structure of their moduli spaces. We intend to construct new invariants that
will distinguish between different connected components of moduli space.
We study invariants that come from pluricanonical embeddings, the corresponding
generic projections to $\CPt$  and branch curves.
More specifically, we investigate the fundamental groups of the complements in
$\CPt$  of those branch curves.

Let $V$  be an algebraic surface of general type pluricanonically embedded in
$\CP^N, \pi\: V \ri \CPt$  be a generic projection, $S_V  (\sbt \CPt)$  be the
corresponding branch curve.
The topological invariants of $\CPt - S_V$  do not change when the complex
structure of $V$ is  changes continuously (for more details see the
introduction of \cite{MoTe6} ).
Thus, one can use such invariants to distinguish the connected components of
the
corresponding moduli space.

	Artin's braid groups $B_n$  are our most important tool.
In recent years these groups are  becoming more and more popular in different
branches of mathematics, such as the theory of Jones polynomials (which is a
current interest of theoretical physicists).
Our former study shows that some nontrivial properties of braid groups are
intimately connected with the existence of nontrivial geometrical objects of
complex geometry and algebraic surfaces.
Thus, our study of algebraic surfaces and their branch curves gives  new
insight in the field of braid groups and vice versa.

The fact that algebraic surfaces are nontrivial geometric objects was
remarkably confirmed  by S. Donaldson who showed that among algebraic surfaces
one can find homeomorphic non-diffeomorphic (simply-connected)  4-manifolds.
In particular, he produced the first counterexamples to the h-cobordism
conjecture in dimension four.
Donaldson's theory was also used  to construct the first examples of
homeomorphic
non-diffeomorphic (simply-connected) algebraic surfaces of general type
(\cite{FMoM}, \cite{Mo4}).

	We expect thatthe connected components of moduli spaces of algebraic surfaces
(of general type) correspond to the principal diffeomorphism classes of
corresponding topological 4-manifolds.
Thus, it is possible that Donaldson's polynomial invariants will distinguish
these connected components.
However, the definition of Donaldson's invariants seems to be too
``transcendental" for direct computation.
We believe that a more direct geometrical approach must be applied.

	An algebraic surface can be considered as a ``Riemann surface" of an algebraic
function of two  variables, that is, as a finite ramified covering of $\CPt.$
For classification problems we can restrict  ourselves to the so-called stable
case, that is, to the case when the corresponding branch curve $S$ in $\CPt$
is
cuspidal (having only nodes and cusps as singularities).
In such an approach, the most important thing is to study the topology of the
complement $\CPt - S_V$, in particular $\pi_1( \CPt - S_V,*)$  and $\pi_ 2(\CPt
- S_V)$  as $\pi_1$-module, etc.

	The first results on such $\pi_1( \CPt - S_V,*)$  were obtained by O. Zariski
in  the thirties.
Let $X_n$  be a non-singular hypersurface in $\CP^3$  of  $\deg n,$  $X_{ab}$
be a projective embedding of $\CP^1 \tm \CP^1$  corresponding to  $a\ell_1 +
b\ell_2,$ where $a,b$  are positive integers  $(\ell_1 = pt \tm \CP', \ell_2 =
pt\tm\CP').$ Let $S_n$  (resp.  $S_{ab}$)  be the branch curve in $\CPt$
corresponding to a generic projection of $X_n$  (resp.  $X_{ab}$)  to $\CPt.$
Zariski proved that $\pi_1(\CPt - S_3,*)\cong \BZ/2 \ast \BZ/3,$  and
$\pi_1(\CPt - S_{1b}) \cong$ braid group of $S^2$  with  $2b$  strings.
	In 1981 B. Moishezon generalized Zariski's result to $X_n,$  proving that
 $\pi_1(\CPt - S_n,*) \cong B_n/$center.

	These results were almost the only ones known about $\pi_1(\CPt - S_V)$'s (in
the stable case)  till the middle of the eighties.
They gave the impression that these groups are very big (in particular, contain
large free groups).
Such an impression was partly responsible for the general belief in the
following
conjecture:
The Galois coverings corresponding to generic projections of
algebraic  surfaces to $\CPt$  have (as a rule) infinite fundamental groups.
This was a partial case of  Bogomolov's conjecture which stated that algebraic
surfaces of general type with positive index  have infinite fundamental
groups. This is equivalent to the following:  for any simply-connected
algebraic
surface of general type the Chern numbers satisfy the inequality $C_1^2 <
2C_2.$

	In 1984 we disproved this conjecture (\cite{MoTe1}), showing that Galois
coverings  $\tX_{ab}$   corresponding to $\Xab \to \CPt$  (generic projection)
have finite fundamental groups.
One can  check that $\pi_1(\tX_{ab})\cong \ker (\pi_1(\CPt - S_{ab}) /\la\G^2
\ra) \overset\psi\to\twoheadrightarrow S_k) $\linebreak $ (\la\G^{2}\ra$
denotes the normal  subgroup generated by squares of ``geometric generators"),
$\psi$  the standard monodromy homomorphism to the symmetric group
$S_k, $\linebreak $ k = \deg \Xab).$   We proved that $\ker \psi$  is a finite
abelian group.   In 1985 we computed $\ker \Psi$  explicitly as an
$S_k$-module.   These results gave the first sign that, in general, the
fundamental groups  $\pi_1(\CPt - S_V,*)$  are not as big and complicated as
the
earlier theorems of Zariski et al. (see above) indicated.

With all the above in mind, one can start to attack the general problem of
explicitly computing $\pi_1(\CPt - S_V,*)$  for many classes of algebraic
surfaces.
Our belief is that in the case of simply-connected surfaces, the kernels
\newline $\ker(\pi_1(\CPt - S_V,*) \overset\psi\to\twoheadrightarrow$ symmetric
group  ($\Psi$ the standard  monodromy homomorphism) are very often solvable
groups.

We start by considering here the branch curve of a Veronese surface.
Completing the computation of$\pi_1$ will allow us to start the investigation
of
$\pi_2(\CPt - S_V)$   (as $\pi_1$-module) which is our next objective.

For $V$  be a projective algebraic surface, $S_V$  the branch curve
corresponding
to a generic  projection of $V$  to  $\CPt.$
$S_V$ is  a cuspidal  curve.
Our method of studying $\CPt - S_V$  is  based on
the explicit formulas for braid monodromy.
In our papers \cite{MoTe5}, \cite{MoTe6} we developed algorithms and methods to
obtain such formulas for many interesting geometrical cases.
These formulas contain nontrivial information which connects some properties of
braid groups to the existence of cuspidal branch curves and the corresponding
algebraic surfaces.
The classical Van Kampen theorem gives an algorithm for deducing a finite
presentation of $\pi_1(\CPt - S_V,*)$  from the explicit knowledge of the braid
monodromy.
At first glance, this presentation looks hopeless.
We manage to work with such presentations, since we have
discovered certain symmetry properties in the braid monodromies (we call them
``invariance properties").
These properties sharply distinguish our subject from its real analog --
classical knot theory, where finite presentations of a knot group usually
cannot
be simplified.

We shall  also use certain new quotients of braid groups and related algebraic
objects described in \cite{MoTe9} which are important ingredients in the
final description of fundamental  groups $\pi_1(\BC^2 - S_V,*),$
as in \cite{MoTe9} and of $\pi_1(\CPt - S_V,*)$ in this paper.
To get the results in \cite{MoTe9} we also used	degenerations of $V_n$ into a
union of simple rational  surfaces which give degenerations of the
corresponding
branch curves into configurations of simple curves ( lines and conics),
(see \cite{MoTe7}).
Such degenerations make it possible to move from the local
analysis of singularities to the global analysis, which is the most difficult
part of the subject, and to compute the braid monodromy related to $S.$
Using the braid monodromy, we computed $\pi_1(\BC^2 - S_V,*)$ in \cite{MoTe9}.

	Our conjecture is that for many classes of simply-connected algebraic surfaces
the fundamental groups $\pi_1(\CPt - S_V)$  will be extensions of symmetric
groups by  solvable groups.
At the present stage of our knowledge it is difficult to predict the precise
answers.
One possibility is that certain general theorems will be proved about
the structure of the solvable groups in question.
Our knowledge of  $\pi_2(\CPt - S_V)$'s is practically non-existent at the
present time.
That means that it is impossible to even formulate  conjectures about them
before our computations  are completed. Still it is possible that certain
general rules define their structure (as $\pi_1$-modules).
The possibilities described above, if realized, will demand  a study of branch
curves of more  complicated
algebraic surfaces, for example, the non-simply-connected ones.  It is known
that many abstract groups (in particular, all finite ones) could be fundamental
groups of algebraic surfaces.

This paper follows our series of papers, Braid Group Techniques, I -- V
(referred to also as BGTI -- BGTV, \cite{MoTe5}, \cite{MoTe6},
\cite{MoTe7}, \cite{MoTe8}, \cite{MoTe9}).
In BGTI we gave an algorithm to compute
the braid monodromy of different line arrangements in $\CPt.$
In BGTII we considered nodal curves and cuspidal branch curves.
In BGTIII we started to treat the Veronese surface separately.
We constructed special degenerations of the Veronese surface to the union of
planes with which we computed in BGTIV the braid monodromy of the branch cuve
of a generic projection of the Veronese surface to $\CPt$.
Using the results in BGTIV we computed in BGTV the fundamental group
$\pi_1(\BC^2-S_V).$

In Section 1 of this paper we discuss the Van Kampen Theorem
which relates the fundamental group of the complement of a curve to the braid
monodromy of the curve.
In Section 2 we prove our main result: a solvability type theorem on
$\pi_1(\BC^2-S_V)$, and in Section 3
 we present the projective  version of the affine results.

	Theorem 2.4 below states that $\pi_1(\CPt-S,*)$ is ``almost'' a solvable group
in the sense that it has a solvable subgroup of finite index.
Such groups are ``small enough'' so that the computation of $\pi_2$ as a model
over $\pi_1$ makes sense.

 \subhead{1. The Braid Group and the Van Kampen
Theorem}\endsubhead

Let $S$ be a cuspidal curve in $\Bbb C^2,$\ $p=\deg S,$\ $\BC_u=\{(u,y)\}.$

The group $\pi_1(\BC-N, u)$ is a free group.

There exists an epimorphism $\pi_1 (\BC_u - S, u_0) \ri\pi_1(\BC^2-S,
u_0),$ so a set of generators for $\pi_1 (\BC_u - S, u_0)$ determines a set
of generators for $\pi_1 (\BC^2-S, u_0).$

There is a classical theorem of Van Kampen from the 30's that all relations
in\linebreak
$\pi_1(\BC^2-S, u_0)$ come from the braid group $B_p$ via the braid monodromy
$\vp_u$ of $S.$ We shall formulate the theorem precisely in 1.2.

We start with the definition of braid group and a half-twist.
\definition{Definition}\ $\underline{\text{Braid group}\ B_n[D,K]}$:\quad Let
$D$
be a closed disc in $\Bbb R^2,$ \ $K\subset D,$ $K$ finite.
Let $B$ be the group of
all diffeomorphisms $\beta$ of $D$ such that $\beta(K) = K\,,\, \beta
|_{\partial D} = \text{Id}_{\partial D}$\,.
For $\beta_1 ,\beta_2\in B$\,, we
say that $\beta_1$ is equivalent to $\beta_2$ if $\beta_1$ and $\beta_2$ induce
the same automorphism of $\pi_1(D-K,u)$\,.
The quotient of $B$ by this
equivalence relation is called the braid group $B_n[D,K]$ ($n= \#K$).
We sometimes denote by $\overline\beta$ the braid represented by $\beta.$
The elements of $B_n[D,K]$ are called braids.\enddefinition
\definition{Definition}\ \underbar{$H(\sigma)$, half-twist defined by
$\sigma$}:\quad  Let
$D,K$ be as above.
Let $a,b\in K\,,\, K_{a,b}=K-a-b$ and $\sigma$ be a simple
path in $D-\partial D$ connecting $a$ with $b$ s.t. $\sigma\cap
K=\{a,b\}.$
Choose a small regular neighborhood $U$ of $\sigma$ and an
orientation preserving diffeomorphism $f:{\Bbb R}^2 \longrightarrow {\Bbb C}^1$
(${\Bbb C}^1$ is taken with the usual ``complex'' orientation) such that
$f(\sigma)=[-1,1]\,,\,$   $ f(U)=\{z\in{\Bbb C}^1 \,|\,|z|<2\}$\,.
Let $\alpha(r),r\geqslant 0$\,,
be a real smooth monotone function such that $
\alpha(r) = 1$ for $r\in [0,\tsize{3\over 2}]$ and
                $\alpha(r) =   0$ for $ r\geqslant 2.$
Let $H(re^{i\theta})=re^{i(\theta+\alpha(r))},$ and let $H(\sigma): D\to D$
be defined by $fHf\1.$\enddefinition
The following lemma is a technical lemma to be used in the simplified Van
Kampen theorem.

\proclaim{Lemma 1.1}
Let $V$ be a half-twist in $B_p[D, K],  u_0 \not\in K.$
Then there exists $A_V, B_V\in\pi_1(D-K,u_0)$ s.t. $A_V, B_V$ can be
extended to a $g$-base of \ $\pi_1(D-K, u_0)$ and $(A_V)V=B_V.$\endproclaim

 \demo{Proof}
\cite{MoTe10}, XIII.1.1 of \cite{MoTe5}. \qed\edm
We use the existence of $A_V, B_V$ in the formulation of the Van Kampen
theorem.

In \cite{MoTe10} and \cite{MoTe5}we also introduced an algorithm for expressing
$A_V,\ B_V $ in terms of $\{\G_j\}.$

To formulate the Van Kampen Theorem we consider the following situations and
use the following notations:

Let $S$ be a cuspidal curve in $\CPt,$ $p=\deg S.$
Let $L$ be a line at infinity,

$\BC^2=\CPt - L.$

Choose coordinates $x, y$ on $\BC^2.$

$\pi:\BC^2\ri\BC$ projection on the first coordinate, $x$-coordinate,
$\BC_x=\pi\1(x).$

$K(x)=\pi^{-1}{(x)}\cap S \qquad (\# K(x)\leq p).$

$N=\{ x|\# K(x) \lvertneqq p\} \ = \ \{c_i\}.$

Let us choose $u\in\BC,$ $u$ real s.t. $x\ll u$ \ $\forall \ x\in N.$

Let $B_p=B[\BC_u, \BC_u\cap S].$

Let $\vp_u: \pi_1 (\BC-N, u)\ri B_p$ the braid monodromy of $S$  with respect
to
$\pi, u.$

The group $\pi_1 (\BC_u - S, u_0)$ is a free group.

$M'=\{x\in S \ | \ \pi |_S$ is not etale at $x\} \quad (\pi(M')=N).$

Assume $\#\BC_x\cap M'=1 \ \ \forall \ x\in N.$
For $x\in N,$ let $x'=\BC_x\cap M',$
The point $x'$ is either a branch point, a node, or a cusp.

Let $\{\delta_i\}$ be a free geometric base ($g$-base) for $\pi_1 (\BC-N, u).$

By Theorem VI.3.3 of \cite{MoTe10} (see also \cite{MoTe9}), for every
$\delta_i$
there exist $V_i$ and $\nu_i$ where $V_i$ is a half-twist and $\nu_i$ is a
number s.t. $\vp_u(\delta_i) =  V_i^{\nu_{i}.}$ Moreover, $\nu_i=1,2,3$ if
$c'_i=$ a branch point, node or a cusp, respectively.

Let $u_0\in\BC_u, u_0\not\in S,$ $u_0$ below real lines far enough s.t.
$B_p$ does not move $u_0.$

$[A, B] = ABA\1 B\1.$

$\la A, B\ra=ABAB\1 A\1 B\1.$

\proclaim{Theorem 1.2 (Van Kampen Theorem)}

Let $S$ be cuspidal curve in $\BC^2.$

Let $ u, u_0, \vp_u, A_V, B_V$ be as above:

 Let $\{\delta_i\}$ be a $g$-base of $\pi_1
(\BC-N,u)$).

Let $\vp_u(\delta_i)=V_i^{(\nu_i)}$\ $V_i$ a half-twist, \ $\nu_i=1,2,3$ (as
above).

Let $\{\G_j\}^p_{j=1}$ be a $g$-base for $\pi_1 (\BC_u -
S, u_0).$
Then $\pi_1=\pi_1(\BC^2 - S, u_0)$ is generated by the images of
$\G_j$ in $\pi_1$ and we get a complete set of relations from those induced
from
$\vp_u(\delta_i)=V^{\nu_i}_i $, as follows:  $A_{V_{i}} = B_{V_{i}}$ when
$\nu_i=1,$
$[A_{V_{i}}, B_{V_{i}}]=1$ when $\nu_i=2,$
$\la A_{V_{i}}, B_{V_{i}}\ra=1$ when $\nu_i=3,$ and $A_{V_{i}}, B_{V_{i}}$ are
expressed in terms of $\{\G_j\}.$\endproclaim

\demo{Proof}\ (See \cite{VK1}, \cite{Z1}).
\enddemo

We shall also quote here an equivalent form of this theorem.

Let $z_{12}=$ the ``shortest" path connecting $q_1$ and $q_2$.
Let $Z=H(z_{12}).$
Since every half-twist of $B_p$ is conjugate to $Z$ (see Proposition VI.4.4,
\cite{MoTe10} or \cite{MoTe5}.), we get $\rho_u(\delta_i) = Q_i\1 Z^{\nu_i}
Q_i$
for some $Q_i\in B_p.$

\proclaim{Theorem 1.3} \ \ \text{\rm(The original version of the Van Kampen
Theorem)}.

Let $S, u, u_0, \{\delta_i\}, Q_i$ be as above.

Let $\{\G_j\}^p_{j=1}$ be a $g$-system of generators for $F=\pi_1 (\BC_u - S,
u_0).$
Then $\pi_1=\pi_1(\BC^2 - S, u_0)$ is finitely presented by the images of
$\G_j$ in $\pi_1$ and the relations which are images of: $Q_i(\G_1)=Q_i(\G_2)$
when $\nu_i=1, Q_i(\G_1) Q_i(\G_2) = Q_i(\G_2) Q_i(\G_1)$ when $\nu_i=2,$ and
$Q_i(\G_1)Q_i(\G_2)Q_i(\G_1) = Q_i(\G_2)Q_i(\G_1)Q_i(\G_2)$ when $\nu_i=3,$ for
each $c_i$ in $N.$\endproclaim

It is easy to see that $((\G_1) H(z_{12}) = \G_2$ (thus for
$Z=H(z_{12},$ we have $\G_1=A_Z,\ \ \G_2=B_Z$). Therefore, $\forall \ i \
(\G_1)Q_i\cdot Q_i\1\cdot Z Q_i = \G_2 Q_i.$ Thus $((\G_1)Q_i) V_i=(\G_2)Q_i.$
So $(\G_1)Q_i=A_{V_{i}}, (\G_2)Q_i=B_{V_{i}}.$
This establishes the equivalence of the different formulations of the theorem.

This theorem is one of the motivations  to compute braid monodromy of cuspidal
curves.
Since the branch curve of a stable generic projection is cuspidal, we can use
this theorem to compute $\pi_1 (\BC^2 - S, u_0)$ whenever $\vp_u$ can be
computed
explicitly.

The next theorem is the projective Van Kampen theorem to be used in
computations of $\pi_1(\CPt-\ov S, *).$
 \proclaim{Theorem 1.4 (Projective Van
Kampen Theorem)} \ Let $\ov S$ be a cuspidal curve in $\CPt,$ transversal to
the
line in infinity. Let $S=\ov S\cap \BC^2.$
Let $\G_1,\dots, \G_p$ be a $g$-base for $\pi_1(\BC_u-S,u_0).$
Then $\pi_1(\CPt-\ov S,*)\simeq
\pi_1(\BC^2-S,*)\bigg/\la\prod\limits_{j=1}^p\G_j\ra.$\ep

In \cite{MoTe5} we proved that for a $g$-base $\{\delta_i\}$ of
$\pi_1(\BC^2-S,u_0)$ we have $\Dl^2_p=\prod\vp_u(\delta_i)$ where $\vp_u$ is
the braid monodromy w.r.t. $S,\pi, u$ and $\Dl^2_p$ is the generator of
$B_p[\BC_u, \BC_u\cap S].$
Moreover, we proved there that the set of such product-forms of $\Dl^2_p$ form
a complete class under Hurwitz equivalence of factorized expressions.
Moreover, from all these equivalent factorizations we try to choose the
product form which will be most useful for fundamental group calculations.
So to calculate the braid monodromy is equivalent to finding certain
factorizations of $\Dl_p^2.$

Invariance properties are results in which we prove that the braid monodromy
factorization of $\Dl^2_{2p}$ is invariant under certain elements of $B_{2p}.$
We look for elements of $B_{2p}$ that will give us equivalent factorizations of
$\Dl_p^2.$
Establishing invariant properties is essential in order to simplify the
calculations which follow from the Van Kampen Theorem.

For defining invariance properties we need the following definitions:

\demo{Definition} \ \underbar{Hurwitz move}

Let $g_1\cdot\dots\cdot g_k= h_1\cdot\dots\cdot h_k$ be two factorized
expressions of the same element in a group $G.$
We say that $g_1\cdot\dots\cdot g_k$ is obtained from $h_1\cdot\dots\cdot
h_k$ by a Hurwitz move if $\exists \ 1\leq p\leq k-1$ s.t. $g_i=h_i, \ i\neq p,
p+1, g_p=h_ph_{p+1}h\1_p$ and $g_{p+1}=h_p$ or $g_p=h_{p+1}$ and $g_{p+1} =
h_{p+1}\1 h_ph_{p+1}.$\enddemo

\demo{Definition} \ \underbar{Hurwitz equivalence of factorized expressions}

Let $g_1\cdot\dots\cdot g_k= h_1\cdot\dots\cdot h_k$ be two factorized
expressions of the same element in a group $G.$
We say that $g_1\cdot\dots\cdot g_k$ is a Hurwitz equivalent to
$h_1\cdot\dots\cdot h_k$ if $h_1\cdot\dots\cdot g_k$
 is obtained from $h_1\cdot\dots\cdot h_k$ by a finite number of Hurwitz moves.
We denote it by $g_1\cdot\dots\cdot g_k \ \underset\text{He}\to\simeq
 \ h_1\cdot\dots\cdot h_k.$\enddemo

\demo{Definition} \ \underbar{A factorized expression invariant under $h$}

Let $g_1\cdot\dots\cdot g_k$ be a factorized expression in $G, h\in G.$
We say that $g=g_1\cdot\dots\cdot g_k$ is invariant under $h$ if
$g_1\cdot\dots\cdot g_k$ is
Hurwitz equivalent to $(g_1)_h\cdot\dots\cdot (g_k)_h,$
 where $(g_i)_n=h\1 g_ih.$\enddemo

Invariance properties are important in view of the following lemma
(\cite{MoTe5}):

\proclaim{Lemma 1.5}\
If a braid monodromy factorization $\Dl^2_p=\prod Z_i$ is invariant
under $h$ then the equivalent factorization $\Dl^2_p=\prod (Z_i)_h$ is also a
braid monodromy factorization.\endproclaim

Since every factor of a braid monodromy factorization induces a relation on
$\pi_1(\BC^2-S)$ by proving invariant propertiesm we get more information
on\linebreak $\pi_1(\BC^2-S),$ (each $(Z_i)_h$ induces a new relation) and thus
it is an essential addition to the Van ~Kampen theorem.
In the next theorem we shall explain how to get new relations from invariance
properties.

\proclaim{Theorem 1.6}
Let $S, u, \rho_u, B_p$ be as above.
Let $\Dl_p=\prod Z_i$ be a braid monodromy factorization of $\Dl^2_p$ w.r.t.
$\vp_u.$
If a subfactorization $\prod\limits^r_{i=s} Z_i$ is invariant under $h,$ and
$\prod\limits^r_{i=s} Z_i$ induces a relation $\G_{i_1} \cdot\dots\cdot
\G_{i_t}$ on $G$ via the Van Kampen method then $(\G_{i_1})_n \cdot\dots\cdot
(\G_{i_t})_n$ is also a relation.
If \ \ $\prod Z_i$ is invariant under $h,$ and if
$R=\G_{i_{1}}\cdot\dots\cdot\G_{i_{t}}$ is a relation on $G=\pi_1(\BC^2-S,
u_0),$ then
$(\G_{i_{1}})h\cdot\dots\cdot(\G_{i_{t}})h$ is also a relation.\endproclaim

\demo{Proof}

Any relation on $\pi_1(\BC^2-S, u_0)$ is a product of the relations induced
by $\prod Z_i,$ via the Van Kampen method.
Thus, it is enough to consider relations induced from the braid monodromy.
Assume that $V^\nu$ is a factor in $\prod\limits^k_{i=2} Z_i.$
Assume that the relation induced by it is:
$$\alignat 2
&A_V=B_V \qquad \qquad &&\nu=1\\
&[A_V,B_V]=1 \qquad \qquad &&\nu=2\\
&\la A_V,B_V\ra=1 \qquad \qquad &&\nu=3.\endalignat$$

Since $\prod\limits^r_{i=s} \ Z_i$ is invariant under $h,$ \
$\prod\limits^r_{i=s} (Z_i)_h$ is also part of a braid monodromy factorization.
In $\prod\limits^r_{i=s} (Z_i)_h$ we have a factor of the form
 $(V^\nu)_h.$
Since $(V^\nu)_h=(V_h)^\nu,$ $(V^\nu)_h$ induces the following relation on
$\pi_1(\BC^2 -S, u):$
$$\alignat 2
&A_{V_{h}}=B_{V_{h}} \qquad \qquad &&\nu=1\\
&[A_{V_{h}}, B_{V_{h}}]=1  &&\nu=2\\
&\la A_{V_{h}}, B_{V_{h}}\ra=1  &&\nu=3.\endalignat$$

Since $A_{V_{h}} = (A_{V})h, \  B_{V_{h}} = (B_{V})h,$\ $(V^\nu)_h$
induces the following relation on $\pi_1(\BC^2-S, u):$
$$\alignat 2
&(A_{V})h=(B_{V})h \qquad \qquad &&\nu=1\\
&[(A_{V})h, (B_{V})h]=1 \qquad \qquad &&\nu=2\\
&\la (A_{V}), (B_{V})h\ra=1 \qquad \qquad &&\nu=3.\endalignat$$

Thus, if $V^\nu$ induces the relation
$\G_{i_{1}}\cdot\dots\cdot\G_{i_{t}}=1$ then $(V^\nu)_h$ induces the
relation $(\G_{i_{1}})h\cdot\dots\cdot(\G_{i_{t}})h=1.$
\qquad   \qed\enddemo
     \subhead 2.
 The topology of affine complements of Veronese branch curves \endsubhead

Let $S_{V_{3}}$ be the branch curve of a generic projection $V_3\ri\CPt$
defined
in Section 2.
For short we shall denote $S_{V_{3}}$ by $S$ in the sequel.

In \cite{MoTe7} we proved a result concerning $G=\pi_1 (\BC^2 - S,*).$
In this section, we shall quote this result and prove further results
concerning $G,$ presenting it as an ``almost''  solvable group.

Let $B_n$ be the braid group, as in Section 1.
In order to formulate these results we need a few definitions.

\demo{Definition}
\underbar{Transversal half-twists}

The half-twists $H(\sigma_1)$ and $H(\sigma_2)$ will be called {\it
transversal} if $C_1$ and $C_2$ intersect transversally in one point which
is not an end point of either of the $\sigma_i$'s.\edm

\demo{Definition}\ $\underline{\tilde B_n,\ \tilde X_i}$

Let $T_n$ be the subgroup of $B_n$ normally generated by $[X,Y]$ for $X,Y$
transversal half-twists.
$\tilde B_n$ is the quotient of $B_n$ modulo $T_n.$
We choose a frame $X_i$ of $B_n.$
We denote their images in $\tilde B_n$ by $\tilde X_i.$\enddemo

We shall use a slightly different presentation for $\tB_9$ than the one induced
from the standard Artin presentation:
$$\tB_9=\la \tT_1,\dots,\tT_9,\ i\ne 4\ra$$
with the following complete list of relations:
$$\alignat 2\\
&[\tT_i,\tT_j]=1\quad && t_i, t_j\ \text{disjoint}\\
&\la\tT_i,\tT_j\ra=1\quad &&\text{otherwise}\\
&[\tT_1,\tT_2\1\tT_3\tT_2]=1\\
&[\tT_5,\tT_8\tT_9\tT_8]=1\endalignat$$
where $T_i$ is a half-twist corresponding to a path $t_i$ and the $t_i$ are
arranged in the following configuration:

\midspace{1.30in}

 \proclaim{Proposition-Definition 2.0}\
$\underline{G_0(n),\tau,u_1}$

Let $A_{n-1}$ be the free abelian group on $w_1,\dots,w_{n-1}.$
Let us define a $\Bbb Z/2$ skew-symmetric form on $A_{n-1}$ as follows:
$$w_i\cdot w_j=\cases 1 \quad & |i-j| =1\\
0 \quad & \text{otherwise}.\endcases$$
There exists a unique central extension $G_0(n),$ of $\Bbb Z/2$ by
$A_{n-1},$ with generators $u_1,\dots, u_{n-1}$ that satisfies
$$\align &1\rightarrow \Bbb Z/2\overset b\to\rightarrow G_0(n)\overset
a\to\rightarrow A_{n-1}\rightarrow 1\\
&a(u_i)=w_i\\
&[u_i,u_j]=b(w_i\cdot w_j)=\cases \tau \quad &|i-j|=1\\
0\quad &\text{otherwise}.\endcases\endalign$$
We always consider $G_0(n)$ with the standard $\tilde B_n$-action as follows:
$$(u_i)_{\tilde X_{k}}=\cases u_i^{-1}\tau&\quad k=i\\
u_ku_i&\quad |i-k|=1\\
u_i&\quad |i-k|\geq 2\endcases.$$
\endproclaim

\proclaim{Claim 2.1}\
$Ab(G_0(n))=A_{n-1}$ (free abelian group on $n-1$ generators),\quad
$G_0(n)'=\{\tau,1\}(\simeq \Bbb Z/2).$\ep

\demo{Proof} Claim III.6.4 of \cite{MoTe9}.\quad
$\square$\edm

Consider the semidirect product, $\tilde B_9\ltimes G_0(9),$ with respect to
the
standard $\tilde B_9$ action on $G_0(9).$

We will work with a more concrete presentation of $G_0(9)$  that will be
compatible with the chosen presentation of $\tB_9:$

 Let $G_0(9)$ be the group
generated by $g_i,$\ $i=1,\dots, 9,$\ $i\neq 4,$ with the following list of
relations.

$[g_1,g_2]^2=1$

$[g_1,g_2]\in\BCen(G_0(9))$

$[g_i,g_j]=\cases 1 \quad & t_i, t_j\ \text{are disjoint}\\
[g_1,g_2]\quad & \text{otherwise.} \endcases$

Denote $\nu=[g_1,g_2].$
Let us reformulate the relations of $G_0(9)$ as follows:

\flushpar $G_0(9) =\big\la g_1,\dots,\check g_4,\dots,
g_9\bigm|[g_i,g_j]=\cases 1\quad & T_i, T_j\ \text{are disjoint}\\
&\qquad\qquad \qquad\qquad\  \tau^2=1\quad \tau_{g_i}=\tau \\
 \tau\  & \text{otherwise}
.\endcases\big\ra$

\demo{Definition}\
$\underline{v_1, N_9, G_9, \tilde\psi_9: G_9\rightarrow S_9}$

$v_1=(\tilde X_2\tilde X_1\tilde X_2^{-1})^2\tilde X_2^{-2}$ for a frame
$X_1,\dots,X_8$ of $B_9.$

$N_9=$ The subgroup normally generated by $c\tau^{-1},\ (u_1v_1^{-1})^3$ \
$(\tau$ an element of $G_0(9),$ (see above) and $c$ an element of $\tilde B_9$
(see above).

$G_9=\tilde B_9\ltimes G_0(9)/{N_9}$

$\tilde \psi_9: G_9\twoheadrightarrow S_9,$ defined by $
\tilde\psi_9(\alpha,\beta)=\tilde\psi_9(\alpha),$ where $\tilde\psi_9:\tB_9\to
S_9$ is the  homomorphism to the symmetric group, induced from the
standard homomorphism $B_9\overset \psi_9\to\rightarrow S_9.$ \enddemo

\proclaim{Proposition 2.2}\ Let $V_3$ be the Veronese surface of order 3.
Let $S_3$ be the branch curve of a generic projection $V_3\ri \CPt.$
Let $\Bbb C^2$ be a big ``affine piece'' of $\CPt.$
Let $S=S_3\cap \CPt.$
 Let $G=\pi_1(\Bbb C^2-S).$
Then $G\cong G_9$ s.t. $\psi: G\rightarrow S_9$ is compatible with
$\tilde\psi_9:G_9\rightarrow S_9.$
\endproclaim
\demo{Proof}
Theorem 6.1 of \cite{MoTe9}.
\edm

Let $\psi_n$ be the standard homomorphisms $B_n\overset
\psi_n\to\rightarrow S_n (=$ symmetric group).
Let $Ab$ be the standard homomorphism $B_n\overset
{Ab}\to \rightarrow\Bbb Z=Ab(\tB_n)=\tB_n/\tB_n'.$
Since $\psi_n([X,Y])=1,$ and $Ab([X,Y])=1,$ $\psi_n$
 and $Ab$ induce  homomorphisms on $\tilde B_n.$
   \demo{Definition}\ $\underline{\tilde\psi_n,\tilde
P_n,\tilde P_{n,0},c}$

$\tilde \psi_n:\tilde B_n\rightarrow S_n,$  the induced homomorphism
from $\psi_n$.

 $\widetilde{Ab}\: \tilde B_n\overset \widetilde{Ab}\to \rightarrow\Bbb Z,$ the
induced homomorphism from $B_n \overset{Ab}\to \rightarrow\Bbb Z.$

 $\tilde P_n=\ker \tilde\psi_n.$

$\tilde P_{n,0}=\ker \tilde\psi_n\cap \ker \widetilde{Ab}= \ker \tilde
{P}_n\ri Ab (\tilde B_n)=\Bbb Z$.

$c=[\tilde X_1^2,\tilde X_2^2]$\quad for 2 consecutive half-twists.
\edm

For the proof of the main result  (Proposition 2.4) we need Proposition 2.3
concerning quotients of the braid group. For a subgroup $H$ denote by $H'$ the
commutatator subgroup of $H.$

\proclaim{Proposition 2.3}\ Let $\tX_i$ be a frame in $B_n.$
Let $c=[\tX_1^2, \tX^2_2].$
Then

$c=[\tX_1^2, \tX_2^2] = [\tX_1^2, \tX^2_{i+1}] = \cdots = [\tX_{n-2}^2,
\tX^2_{n-1}].$

Moreover, $(\tP_n)' = (\tP_{n,0})' = \{1, c\} \simeq \BZ_2.$

$Ab(\tP_n) =$ free abelian group on $n$ generators;
$\xi_1,\dots,\xi_{n-1},\tX_1^2,$ where $\xi_1=(\tX_2\tX_1\tX_2\1)^2\tX_2^{-2}$
and $\xi_i$ is conjugate to $\xi_1,$\ $\xi_i\in\tP21
_{n,0}.$

$\tB_n$ acts on $\tP_{n,0}$ by conjugation.

$\tP_{n,0}$ with this action is isomorphic to $G_0(n)$ with the standard
$\tB_n$-action as defined previously.

There exists a series: $1\sbq (\tP_{n,0})' \sbq \tP_{n,0} \sbq\tP_n\sbt \tB_n$
s.t. $\tB_n / \tP_n = S_n,$\linebreak $ \tP_n/\tP_{n,0}\simeq \BZ,$\quad $
\tP_{n,0}/(\tP_{n,0})' \simeq A_{n-1} \simeq Ab (G_0(n)),\ \ (\tP_{n,0})'\simeq
\BZ_2.$\ep

\demo{Proof} Theorem III.6.4 of \cite{MoTe9}.\ See  \cite{MoTe5}, Chapters 4,
5 for more information about $\tP_n$ and $\tP_{n,0}$. \ \qed
\edm

The following theorem is the main result of this paper.
It is a structure theorem for
$G=G_9,$ which states that $G_9$ is an almost solvable group.
Using this result, one can start to compute $\pi_2$ as a model over $\pi_1.$

\demo{Definition} \ $\underline{Ab_9,\ H_9, H_{9,0}}$

$Ab_9=$ Abelization map of $G_9.$

$H_9=\ker\tilde\psi_9$

$H_{9,0}=\ker\tilde\psi_9\cap \ker Ab_9$\edm
\proclaim{Proposition 2.4}\
There exists a series
$1\subseteq H_{9,0}'\subset H_{9,0}\subset H_9\subset G_9,$
where $G_9/H_9\simeq S_9,$ \quad\qquad $H_9/H_{9,0}\simeq\Bbb Z,$\quad
\qquad  $H_{9,0}/H'_{9,0}\simeq
(\Bbb Z+\Bbb Z/3)^8,$\linebreak $H_{9,0}'=H_9'=\{1,c\}\cong \Bbb Z/2,$ where
$c\in\BCen(G_9).$
Moreover, $Ab(G_9)=\BZ.$
\endproclaim
\demo{Proof}
Let $\tilde T_i$\ $i=1\dots 9\ i\ne  4$ be the base of $\tB_9$ as above.
Let $g_i,$\ $i=1,\dots,9,$\ $i\ne 4,$ be the chosen base
for $G_0(9).$

 Let

$c=[\tT_1^2,\tT_2^2]$

$\tau=[g_1,g_2].$

$\tB_9$ acts on $G_0(9)$ as follows:
$$
(g_i)_{T_{k}}=\cases g_i\quad & i,k\ \text{disjoint}\\
g_i^{-1}\tau\quad & i=k\\
 g_ig_k^{- 1}\ \text{or}\ g_kg_i\quad &
\text{otherwise}.
\endcases$$

It is easy to check that we have the following relations in $\tB_9\ltimes
G_0(9)$

$$[\tT_i^2,\tT_j^2]=\cases 1\quad & i,j\ \text{disjoint or}\ i=j\\
c\quad &  \text{otherwise}\endcases$$
$$[g_i,g_j]\ =\ \cases 1\quad & i,j\ \text{disjoint or}\ i=j\\
\tau\quad &  \text{otherwise}\endcases$$
$$\align &c_j\tau\in\BCen \tB_9\ltimes G_0(9)\qquad\qquad\qquad\\
& G_0(9)'=\{1,\tau\}\qquad\qquad\\
&c^2=\tau^2=1\qquad\qquad\qquad\\
&c\ne\tau.\endalign$$
We recall from the proof of Proposition 2.3 (which appears in \cite{MoTe9})
that
$\tP_{9,0}$ is generated by $\xi_1,\dots,\xi_9\quad i\ne 4,$\quad $\tP_9$ is
generated by $\tP_{9,0}$ and $\tT_1^2,$\quad $\tP_9'=\tP_{9,0}'=\{1,c\},$
where
$$[\xi_i,\xi_j]=\cases 1\quad & i,j\ \text{disjoint},\ i=j\\
c\quad & \text{otherwise}
\endcases$$
$(\xi_1=(\tT_2\tT_1\tT_2\1)^2\tT_2^{-2}$ and $\xi_i$ is conjugate to $\xi_1).$

Let $\zeta_i=g_i\xi_i\1.$

Since $\tT_1,\tT_2$ can be extended to a frame of $\tB_9,$
$\xi_1$ can be considered as $v_1$ in the above definition.
{}From the above commutators one can see that $\zeta_i^3$ is
conjugated to $\zeta_1^3$, up to multiplication by $c\tau.$
Thus $N_9$ can be
represented as $$N_9=\la c\tau, \zeta_i^3\ i=1\dots 9,\ i\ne 4\ra.$$

To prove the theorem we consider first another quotient of $\tB_9\ltimes
G_0(9):$

Let
$$C_9=\la c\tau\ra$$
$$\hat G_9=\tB_9\ltimes G_0(9)/ C_9.$$
There exist natural epimorphisms $\hat I,\ \hat J.$
$$\tB_9\ltimes G_0(9)\overset{\hat I}\to \twoheadrightarrow \hat G_9,\quad \hat
G\overset{\hat J}\to\twoheadrightarrow G_9,$$
s.t. $\hat I, \ \hat J$ in the natural epimorphism $\tB_9\ltimes G_0(9)\to
\tB_9\ltimes G_0(9)/N_9=G_9.$ Let $$M_9=\la \zeta_i^3\quad i=1,\dots,9\quad
i\ne
4\ra.$$ Let $\hat M_9$ be the image under $\hat I$ of $M_9$ in $\hat
G_9.$\enddemo \demo{Remark}
By abuse of notation we use the same notation for elements of $\tB_9,\ G_0(9)$
and their images in $\tB_9\ltimes G_0(9)$ and its different quotients.\enddemo
We divide the proof into 11 claims and corollaries.
\demo{Claim 1}

(a)\ $\ker \hat I=C_9\simeq \BZ_2.$

(b)\ $\hat I(N_9)=\hat M_9,$\quad $N_9/C_9\simeq \hat M_9.$

(c)\ $\hat M_9$ is a normal commutative subgroup of $\hat G_9.$

(d)\ The induced homomorphism $G_9=\tB_9\ltimes G_0(9)/N_9\overset
I\to\rightarrow \hat G_9/\hat M_9$ is an isomorphism.

(e)\ $\hat M_9=\ker(\hat G_9\overset \hat J\to\rightarrow G_9),$\quad $\hat
G_9/\hat M_9=G_9.$

(f)\ There exists a short exact sequence
$$1\to \hat M_9\to \hat G_9\overset \hat J\to\rightarrow G_9\to 1.$$
\edm
\demo{Proof}

(a) \ Clearly $C_9=\ker\hat I.$ Since $c,\tau\in\BCen \tB_9\ltimes
G_0(9)$ and $ c^2=\tau^2=1$ then $(c\tau)^2=1$ since $c\ne \tau,\quad
C_9=\BZ_2.$

(b)\ $N_9=\la M_9, C_9\ra\Rightarrow I(N_9)=\la \hat I(N_9),\hat I(C_9)\ra=
(\hat M_9,1\ra=\hat M_9.$
Since $\ker \hat I=C_9\subseteq N_9,\quad \ker\left(\hat
I\bigm|_{N_{9}}\right)=C_9.$
Thus $N_9/C_{9}=\hat M_9.$

(c)\ Consider $\tB_9\ltimes G_0(9).$
By Lemma IV.6.1 of \cite{MoTe9}:
$$[\tT_j^{\pm 2},g_i^{\pm 1}] =\cases 1\quad & t_i,t_j\ \text{disjoint},\ i=j\\
c\quad & \text{otherwise}.\endcases$$
By Claim II.4(e) of \cite{MoTe9}
$$[T_j^{-2}T_k^2,g_i^{\pm 1}]=\cases 1\quad & t_i,t_k\ \text{disjoint}\\
c\quad & \text{otherwise}.\endcases$$
By Corollary IV.4.2 for every $\xi_\ell$ and every $i,$ s.t. $t_i\cap
t_\ell\ne\emptyset \ \exists k$ s.t. $t_k\cap t_i\ne\emptyset$ and
$\xi_\ell=T_j^{\pm 2}T_k^{\pm 2}.$
Thus
$$[\xi_\ell^{\pm 1},g_i^{\pm 1}]=\cases 1\quad & i,\ell\ \text{disjoint}\\
c\quad & \text{otherwise}.\endcases$$
On the other hand,
$$[g_\ell^{\pm},g_i^{\pm}]=\cases 1 \quad & t_i, t_\ell\ \text{disjoint}\\
\tau\quad & \text{otherwise}.\endcases$$
Since $c,\tau$ are central elements and $\zeta_\ell=g_\ell\xi_\ell\1,$ the last
2 equations imply: $$[\zeta_\ell^{\pm 1},g_i^{\pm 1}]=\cases 1\quad &
t_i,t_\ell\
\text{disjoint}\\ c\tau\quad & \text{otherwise}.\endcases$$
So $$[\zeta_\ell^{\pm 3},g_i^{\pm 1}]=\cases 1\quad & t_i,t_\ell\
\text{disjoint}\\
c\tau\quad & \text{otherwise}.\endcases$$
Therefore, $\zeta_\ell^3$ and $g_i$ commute in $\tB_9\ltimes G_0(9)/C_9.$
The same arguments imply that $\zeta_\ell^3$ and $\xi_i$ commute.
Thus, $\zeta_\ell^3$ commutes with $\zeta_i,$ and therefore with $\zeta_i^3.$
Thus, $\hat M_9$ is commutative

(d)\ $\ker(\tB_1\ltimes G_0(9)\twoheadrightarrow \hat G_9)=C_9.$
Thus $\ker\left(\tB_9\ltimes G_0(9)\twoheadrightarrow \hat G_9/\hat M_9\right)=
\la C_9,M_9\ra=N_9.$
By the isomorphism theorem we get (d).

(e)\  $\ker(\tB_9\ltimes G_0(9)\to \tB_9\ltimes G_0(9)/N_9)=N_9.$
Thus $\ker(\tB_9\ltimes G_0(9)/C_9\to \tB_9\ltimes G_0(9)/N_9)=N_9/C_9$ which
equals $\hat M_9.$

 (f) From (e).
\quad  \qed\edm

Thus we have
$$\alignat4
&\tB_9\rtimes G_0(9)&&\ \overset{\hat I}\to\rightarrow&&\ \hat G_9&&\
\overset\hat J \to G_9\\ &\quad \cup&&  && \ \cup\\
&\quad N_9 &&\ \to\ &&M_9&&\ \to 1\\
&\tB\rtimes G_0(9)/N_9 &&\ \overset I\to{\osr } &&\ \hat G_9/\hat M_9 &&\
\overset J\to\rightarrow G_9\endalignat$$

\demo{Claim 2}

$Ab(G_9)=Ab(\hat G_9)=Ab(\tB_9\rtimes G_0(9))=Ab(\tB_9)=Ab(B_9)=\BZ.$
\edm
\demo{Proof}
Recall that conjugate elements are equal under abelization $(Ab(G)=G/G').$
Since $B_n$ is generated by half-twists and they are all conjugate to each
other, $Ab(B_n)$ is generated by one element, whose infinite order remains
valid under abelization.
So $Ab(B_9)=\BZ.$
Since $\tB_9=B_9/\text{subgroup of}\ B_9',$ \ $Ab(\tB_9)$ is also $\BZ.$
{}From the action of $\tB_9$ on $G_0(9)$ one can see that for every $i,$\ $i\le
9,\
i\ne 4,$ $\exists k,$ \ $i\ne k,$ s.t. $(g_i)_{\tilde T_{k}}=g_ig_k^{\pm 1}.$
Under abelization this relation becomes $g_i=g_ig_k^{\pm}.$
Therefore, $g_k^{\pm 1}=1.$
So in the abelization of $\tB_9\ltimes G_0(9),$ the elements of $G_0(9)$ are 1.
Thus $Ab(\tB_9\ltimes G_0(9))=\BZ.$
Since $C_9\subseteq(\tB_9\ltimes G_0(9))',$\ $Ab(\tB_9\ltimes G_0(9)/C_9)$ is
also $\BZ.$
So $Ab(\hat G_9)=\BZ.$
Since $M_9$ consists of degree $0$ elements, $Ab(\hat G_9/\hat M_9)=\BZ.$
Thus $Ab(G_9)=\BZ.$ \qed
\edm
Let $\hat H_9: \ker \hat G_9\overset\hat\psi_9\to\twoheadrightarrow S_9 (=$ the
symmetric group on 9 elements).
\demo{Claim 3}

(a) $\hat H_9 \simeq  \tilde P_9\ltimes G_0(9)/C_9.$

(b) There exists an epimorphism $\hat H_9\twoheadrightarrow G_0(18)$ which is
compatible with $Ab(\hat H_9)\to A_{17}$\ $(A_{17}=Ab\ G_0(18)).$

(c) $Ab(\hat H_9)$ is freely generated by $\xi_1\dots \xi_9$\ $i\ne 4,$ \
$g_1\dots  g_9$ \ $i\ne 4,$ $\tilde T_1^2.$

(d) $M_9\subset \BCen \hat H_9.$

(e) $H_9\simeq\tilde P_9\ltimes G_0(9)/N_9.$

(f) $\hat H_9$ maps onto $H_9$ under $\hat G_9\overset{\hat J}\to\rightarrow
G_9.$

(g) $\hat H_9/\hat M_9\simeq H_9.$
\edm \demo{Proof}

(a) $\ker(\tB_9\to S_9)=\tilde P_9,$ and $\tB_9\ltimes G_0(9)\to S_9$ factors
through $\tilde B_9.$
So\linebreak $\ker(\tB_9\ltimes G_0(9)\to S_9)$ is $\tP_9\ltimes G_0(9).$
So $\hat H_9=\ker(\tB_9\ltimes G_0(9)/C_9\to S_9)=\tilde P_9\ltimes
G_0(9)/C_9.$

(b) (c) We shall skip the precise proof of (b) and (c).
The idea is to derive the presentation of $\hat H_9$ from a
presentation of $\tP_9,$ (Proposition 2.3) and Proposition-Definition 2.0.

(d) The relevant relations are to be found in the proof of Claim 3(c).

(e) Similar to (a).

(f) (g) $\hat G_9\overset\hat\psi_9\to\rightarrow S_9$ factors through
$G_9\overset\psi_9\to\rightarrow S_9,$ i.e., $ \oversetbrace \hat\psi_9\to
{\hat G_9\rightarrow G_9\overset \psi_9\to\rightarrow S_9}.$
So $\ker\hat\psi_9$ is mapped into $\ker \psi_9.$
We can see it differently using (a) and (e):
 Since $C_9\subset N_9,$ there exists a natural epimorphism (restriction of
$\hat J)$
$$\tP_9\ltimes G_0(9)/C_9\twoheadrightarrow \tP_9\ltimes G_0(9)/N_9.$$
By (a) and (c) this map is actually
$$\hat H_9\twoheadrightarrow  H_9.$$
The kernel of this map is $N_9/C_9,$ which is $\hat M_9$ by Claim 1(b).
\quad \qed

\demo{Claim 4}

(a) $\hat M_9$ is isomorphic to its image in $Ab(\hat H_9).$

(b) $\hat M_9\cap\hat H_9'=\{1\}.$\edm
\demo{Proof}
By Claim 3 (c), $Ab(\hat H_9)$ is freely generated by 17 elements
$\xi_1,\dots,\xi_9,$\ $i\ne 4,$ \ $g_1,\dots,g_9,$\ $i\ne 4,$ \ $\tilde T_1^2,$
 regarded as elements of the abelization.

$M_9$ is generated by $\zeta_1,\dots,\zeta_9,$\ $i\ne 4$ and
$\zeta_i=\xi_ig_i\1.$
The image of $\hat M_9$ in the abelization is generated by
$\{g_i\xi_i\1\}_{i=1\ i\ne 4}^9$ which are 8 products of different elements of
a free base.
Thus $M\simeq $ Image of $M$ in $Ab(\hat H_9).$
Thus $\ker Ab_{\hat H_{9}}\bigm|_{\hat M_{9}}=1.$
But this kernel is $\hat M_9\cap \hat H_9'.$ \qed\edm
\demo{Claim 5}\ $\hat H_9'\simeq H_9'.$
\edm
\demo{Proof}
Consider $\hat G_9\twoheadrightarrow G_9$ with $\ker \hat M_9.$
Under this map $\hat H_9\to H_9$ (Claim 3(f)) and $\hat H_9'\to H_9'.$
Evidently, $\ker\hat H_9'\to H_9'$ is $\hat H_9'\cap \hat M_9$ which is $\{1\}$
by the previous claim. \qed\edm

\demo{Claim 6}\ $(\hat H_9')\simeq\BZ_2.$
\edm
\demo{Proof} By Claim 3(a), $\hat H_9=\tilde P_9\rtimes G_0(9)/C_9.$
$\tilde P_9'=\{1,c\},$\ $(G_0(9))'=\{1,\tau\}.$
Thus $(\tilde P_9\rtimes G_0(9))'=\{1,c,\tau\}.$
When dividing by $c=\tau$ \ $(C_9=\la c\tau\ra)$ we get $(\hat
H_9)'=\{1,\tau\}\subseteq \BZ_2.$
We have to prove that it does not collapse completely.
We have:
$$1\to(\hat H_9)'\to \hat H_9\to Ab(\hat H_9)\to 1.$$
Consider the short exact sequence from Proposition-Definition 2.0:
$$1\to\BZ_2\to G_0(18)\to A_{17}\to 1.$$
By Claim 3(c) we have $Ab(\hat H_9)\simeq A_{17}.$
By Claim 3(b) there exists $\hat H_9\twoheadrightarrow G_0(18).$
Thus there exists $(\hat H_9)'\twoheadrightarrow \BZ_2.$
So $\hat H_9'\simeq \BZ_2.$ \qed\edm
\proclaim{Corollary 7}

$H_9'=\BZ_2.$\ep
Similar to Claims 3, 5, 6, we have:
\demo{Claim 8}

(a) $Ab(\hat H_{9,0})$ is freely generated by $g_1,\dots g_9,$ \ $i\ne 4,$\
$\xi_1,\dots,\xi_9,$\ $i\ne 4.$

(b) $\hat H_{9,0}/\hat M_9\simeq H_{9,0}.$

(c) $(\hat H_{9,0})'=H_{9,0}'.$

(d) $(\hat H_{9,0})'=\BZ_2.$
\edm
\proclaim{Corollary 9}

$H_{9,0}'=H_9'\simeq\BZ_2.$\ep
\demo{Claim 10}

$H_{9,0}/H_{9,0}'=(\BZ_3\oplus\BZ)^8.$\edm
\demo{Proof}
By Claim 8(c), $\hat H_{9,0}/\hat M_9\simeq H_{9,0}.$
Thus $H_{9,0}/H_{9,0}'=Ab(H_{9,0})=$\linebreak $Ab(\hat H_{9,0}/\hat M_9).$
Since $\hat M_9$ is isomorphic to its image in $Ab(\hat H_9)$ (Claim 4(a)), we
get\linebreak $Ab(\hat H_{9,0})/\hat M_9.$
We can take $\{g_i,\zeta_i\}_{i=1\ i\ne 4}^9$ as free generators for $Ab(\hat
H_{9,0}).$
Thus $Ab(\hat H_{9,0})/\hat M_9=\left(\sum\limits^g\Sb i=1\\ i\ne 4\endSb
(g_i)\oplus \sum\limits^g\Sb i=1\\ i\ne
4\endSb\la\zeta_i\ra\right)\bigg/\sum\limits^9\Sb i=1\\i\ne
4\endSb\la\zeta_i^3\ra= \sum\limits\Sb i=1\\ i\ne 4\endSb\la g_i\ra
\oplus\sum\limits^9\Sb i=1\\ i\ne 4\endSb\left( \la \zeta_i\ra\big/
\la\zeta_i^3\ra\right)$\linebreak $\simeq(\BZ\oplus\BZ_3)^8.$ \qed\edm
 \proclaim{Corollary 11}

$H_9\big/ H_{9,0}=\BZ.$\endproclaim
\demo{Proof}
$H_{9,0}=\ker(H_9\to Ab(\hat G_9)).$
Since $Ab\ \hat G_9\simeq Ab\ H_9=\BZ,$ we have
$H_{9,0}=\ker(H_9\twoheadrightarrow \BZ).$\qed\edm
The proof of the different statements of the theorem are Claims and Corollaries
11, 10, 9, 2, and the definition of $H_9$ for $G_9/H_9\simeq S_9.$\qed\ \ \
for Theorem 2.4\edm

\subhead{3.  The fundamental group of complements of a Veronese branch
curve in $\CPt$}  \endsubhead

\proclaim{Theorem 3.1}\ Let $V_3$ be the Veronese surface of order 3.
Let $\ov S$ be the branch curve of a generic projection $V_3\ri \CPt.$
Let $\ov{G} = \pi_1(\CPt - \ov S).$
Then there exist $w_0\in H_{9,0}$ s.t. $\ov{G}\simeq \ov{G}_9 = G_{9}/\la
X_1^{18} w_0\ra.$\ep

\proclaim{Theorem 3.2}\   Let $\overline H_9$ and $\overline H_{9,0}$ be the
images of $H_{9,0}$ and $H_q$ in $\overline G_9.$
Then $\overline H_{9,0}'=\overline H_q'$ and
$$1\subseteq \overline H_{9,0}'\subseteq \overline H_{9,0}\subseteq
\overline H_q\subset \overline G_9$$
where
$$\overline G_9/\overline H_9\simeq S_9,
\quad
\overline H_9/\overline H_{9,0}\simeq \Bbb Z_9,\quad
\overline H_{9,0}/\overline H_{9,0}'\simeq (\Bbb Z+\Bbb Z/3)^8,\quad
\overline H_{9,0}'=\ov H_9'\simeq \Bbb Z/2.$$\endproclaim

\end

{}.

\end